\documentclass[a4,12pt,leqno]{article}
\usepackage{amsmath,amssymb,amsfonts}
\usepackage[english]{babel}
\usepackage{epic,eepic}
\usepackage{amscd,epsfig,amsthm,graphicx}
\usepackage[flushleft]{paralist} 
\usepackage{ascmac} 
\usepackage[all]{xy} 
\usepackage{array} 
\makeindex			 
\newtheorem{theorem}{Theorem}
\newtheorem{lemma}[theorem]{Lemma}

\theoremstyle{definition}
\newtheorem{definition} 
			{Definition}

\newtheorem{remark}{Remark}

\makeatletter
    
    \@addtoreset{equation}{section}
\makeatother

\def\N{{\mathbb N}}
\def\Z{{\mathbb Z}} 
\def\C{{\mathbb C}} \def\R{{\mathbb R}}

\def\cL{{\mathcal L}}

\begin{document}
\title{ {\Huge Analytic properties of generalized Mordell-Tornheim 
		type of multiple zeta-functions and $ L $-functions}}


\author{Takashi Miyagawa
\footnote{Graduate School of Mathematics, Nagoya University, 
Chikusa-ku, Nagoya 464-8602, Japan
 (Email:d15001n@math.nagoya-u.co.jp).}}

\date{}


\maketitle

\begin{abstract}
 Analytic properties of three types of multiple zeta functions, that is, the
 Euler-Zagier type, the Mordell-Tornheim type and the Apostol-Vu type
 have been studied by a lot of authors.
 In particular, in the study of multiple zeta functions of the
 Apostol-Vu type, a generalized multiple zeta function, 
 including both the Euler-Zagier type and the Apostol-Vu type, 
 was introduced.
 In this paper, similarly we consider generalized multiple 
 zeta-functions and $ L $-functions, which include both the 
 Euler-Zagier type and the Mordell-Tornheim type as special cases. 
 We prove the meromorphic 
 continuation to the multi-dimensional complex space, and give 
 the results on possible singularities.
\end{abstract}

\section{Introduction}
The Euler-Zagier type of multiple zeta-function $ \zeta_{EZ,r} $ is
defined by
\begin{eqnarray}\label{EZ-zeta}
       	\zeta_{EZ,r}(s_1,\cdots,s_r)
               &=& \sum_{1 \leq m_1 < \cdots < m_r}
               	\frac{1}{m_1^{s_1} m_2^{s_2} \cdots m_r^{s_r}}	\nonumber\\
               & =& \sum_{m_1=1}^\infty \cdots \sum_{m_r=1}^\infty
                       \frac{1}{m_1^{s_1} (m_1 + m_2)^{s_2} \cdots 
                       (m_1 + \cdots + m_r)^{s_r}},
\end{eqnarray}
where $ s_1,s_2,\cdots,s_r $ are complex variables, 
and the series (\ref{EZ-zeta}) is absolutely convergent in the region
\begin{equation*}
       	\{ (s_1,\cdots,s_r) \in \C^r \mid
               \mathrm{Re} (s_{r-k+1} + s_{r-k+2} + \cdots + s_r) > k 
               \ (k = 1,2. \cdots ,r) \}.
\end{equation*}
The Mordell-Tornheim type and Apostol-Vu type of multiple zeta-functions
are defined by
\begin{equation}\label{MT-zeta}
       	\zeta_{MT,r}(s_1,\cdots,s_r;s_{r+1})
            =  \sum_{m_1=1}^\infty \cdots \sum_{m_r=1}^\infty
                \frac{1}{m_1^{s_1} \cdots m_r^{s_r}
                (m_1 + \cdots + m_r)^{s_{r+1}}}
\end{equation}
and
\begin{equation}\label{AV-zeta}
	\zeta_{AV,r}(s_1,\cdots,s_r;s_{r+1})
              = \sum_{1 \leq m_1 < \cdots < m_r }
                \frac{1}{m_1^{s_1} \cdots m_r^{s_r}
                (m_1 + \cdots + m_r)^{s_{r+1}}}
\end{equation}
where $ s_1,\cdots,s_r,s_{r+1} $ are complex variables. The series
(\ref{MT-zeta}) and (\ref{AV-zeta}) are absolutely convergent in 
\begin{equation}\label{eq:abs_MT_AV}
       	\{ (s_1,\cdots,s_r,s_{r+1}) \in \C^{r+1} \mid
        \mathrm{Re} (s_j) > 1 \ (1 \leq j \leq r), \ \ 
        \mathrm{Re}(s_{r+1}) > 0
	\}.
\end{equation}

For the meromorphic continuation to the whole space $ \C^r $ of (\ref{EZ-zeta}), 
Akiyama Egami and Tanigawa \cite{AET} and Zhao \cite{Zhao}, proved independently 
of each other. Matsumoto \cite{Mat} gave an alternative proof of the analytic 
continuation using the Mellin-Barnes integral formula
\begin{equation}\label{eq:Mellin-Barnes}
       	(1 + \lambda )^{-s} = \frac{1}{2\pi i} 
        		\int_{(c)} 
			\frac{\Gamma (s + z) \Gamma(-z) }
                        {\Gamma(s)} \lambda^{z}  dz,
\end{equation}
where $ s, \lambda \in \C , |\arg \lambda| < \pi , \lambda \neq 0 $, 
and $ c \in \R $ , $ -\mathrm{Re}(s) < c < 0 $ and
the path of integration is the vertical line from $ c - i \infty $ to
$ c + i \infty $. Also, Matsumoto \cite{Mat2} proved the meromorphic 
continuation in the same way for (\ref{MT-zeta}) and (\ref{AV-zeta}). 
In particular, Matsumoto introduced the following function
in the process of proving the meromorphic continuation of (\ref{AV-zeta}).
Let $ 1 \leq j \leq r $, and define
\begin{eqnarray}\label{GAV-zeta}
  	\lefteqn{\widehat{\zeta}_{AV,j,r}(s_1,\cdots,s_j; s_{j+1},
        	\cdots, s_r; s_{r+1})} \nonumber\\
  	&& \qquad \qquad = \sum_{1 \leq m_1 < \cdots < m_r}
                \frac{1}{m_1^{s_1} \cdots m_r^{s_r}
                (m_1 + \cdots + m_j)^{s_{r+1}}} 
                \label{eq:Matsumoto2}
\end{eqnarray}
where $ s_1,\cdots,s_r,s_{r+1} $ are complex variables. 
Since $ \widehat{\zeta}_{AV,r,r} = \zeta_{AV,r} $ and
\begin{eqnarray*}
     \widehat{\zeta}_{AV,1,r}(s_1 ; s_2,\cdots, s_{r+1})
        = \zeta_{EZ,r}(s_1 + s_{r+1}, s_3,\cdots ,s_r),
\end{eqnarray*}
(\ref{GAV-zeta}) forms a generalized class including as special cases both the
Euler-Zagier type (\ref{EZ-zeta}) and the Apostol-Vu type (\ref{AV-zeta}).
He, through the recursive structure
\begin{equation}\label{AV-RS}
		\zeta_{AV,r} \ = \  \widehat{\zeta}_{AV,r,r} \rightarrow \ 
        		\widehat{\zeta}_{AV,r-1,r} \rightarrow \ 
                        \widehat{\zeta}_{AV,r-2,r} \rightarrow \ 
                        \cdots \ \rightarrow \ 
                        \widehat{\zeta}_{AV,1,r} 
                        \ = \  
                        \zeta_{EZ,r} 
\end{equation}
(here $ A \rightarrow B $ means that $ A $ can be expressed as an integral
involving $ B $; see (\ref{eq:induction}), (\ref{eq:induction1}) and 
(\ref{eq:induction2}) below), discussed analytic properties of those functions.

As an analogue of (\ref{GAV-zeta}), in this paper we define the following 
function, and prove the results on meromorphic continuation and
singularities. The results will be stated in Section 2.

\begin{definition}
Let $ 1 \leq j \leq r $, and define
  \begin{eqnarray}
	&&\widehat{\zeta}_{MT,j,r}(s_1,\cdots,s_j;s_{j+1},\cdots,s_{r+1})
        \nonumber \\
        &&\qquad =  \sum_{m_1=1}^\infty \cdots \sum_{m_r=1}^\infty
        	\frac{1}{m_1^{s_1} \cdots m_j^{s_j}
                (m_1 + \cdots + m_j)^{s_{j+1}} 
                \cdots (m_1+ \cdots +m_r)^{s_{r+1}}},
     \label{eq:Miyagawa1} 
  \end{eqnarray}
where $ s_1, \cdots,s_r, s_{r+1} $ are complex variables.
\end{definition}
Since $ \widehat{\zeta}_{MT,r,r} = \zeta_{MT,r} $ and
\[
        \widehat{\zeta}_{MT,1,r}(s_1,\cdots,s_j;s_{j+1},\cdots,s_{r+1})
        = \zeta_{EZ,r}(s_1+s_2, s_3, \cdots , s_{r+1}),
\]
we see that (\ref{eq:Miyagawa1}) forms a generalized class including as specal cases 
both the Euler-Zagier type (\ref{EZ-zeta}) and the Mordell-Tornheim type 
(\ref{MT-zeta}), which can be illustrated in the following figure.

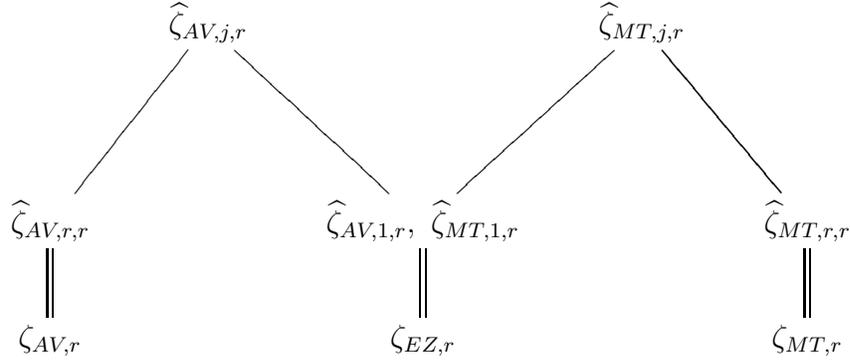
\begin{equation*}\label{eq:xy-fig}
\xymatrix
{
& & \widehat{\zeta}_{AV,j,r} \ar@{-}[ddr] & &
 \widehat{\zeta}_{MT,j,r} \ar@{-}[ddr] &		\\
&&&&&							\\
& \widehat{\zeta}_{AV,r,r}\ar@{=}[d] \ar@{-}[uur] & 
& \widehat{\zeta}_{AV,1,r} ,\  \widehat{\zeta}_{MT,1,r} \ar@{=}[d]
\ar@{-}[uur]
& &  \widehat{\zeta}_{MT,r,r} \ar@{=}[d] \ar@{-}[uul]	\\
& \zeta_{AV,r} & & \zeta_{EZ,r} & & \zeta_{MT,r}	\\
}
\end{equation*}
The series (\ref{eq:Miyagawa1}) is absolutely convergent in the region
\begin{eqnarray*}
	R_{j,r}=
        \left\{
        (s_1,\cdots,s_r,s_{r+1}) \in \C^{r+1} \left|
        \begin{array}{c}
        \mathrm{Re}(s_{r+2-k}+s_{r+3-k}+ \cdots + s_{r+1}) > k \\
        \qquad\qquad\qquad\qquad  (k = 1,2,\cdots,r-j)	\\
        \mathrm{Re}(s_{j+1} + s_{j+2} + \cdots + s_{r+1}) > r-j \\
        \mathrm{Re}(s_{\ell}) > 1 \ \ (\ell = 1,2, \cdots, j)
        \end{array}
        \right.
        \right\},
\end{eqnarray*}
therefore $ \widehat{\zeta}_{MT,j,r} $ is a regular function in $ R_{j,r} $.
This fact can be proved by the evaluation 
\[
	\sum_{m=1}^\infty \frac{1}{(m+N)^\sigma}
        < \int_0^\infty \frac{dx}{(x+N)^\sigma} 
        = \frac{1}{\sigma - 1} \frac{1}{N^{\sigma-1}} \ \ (\sigma > 1)
\]
and the result on the absolutely convergent region (\ref{eq:abs_MT_AV}).

Furthermore, we introduce the following $ L $-function which is a
$ \chi- $analogue of (\ref{eq:Miyagawa1}), and we obtain the results on 
meromorphic continuation and singularities. The results will be stated 
in Section 2.

\begin{definition}\label{def:L_MT,j,r}
Let $ \chi_1, \chi_2 \cdots, \chi_r $ be Diriclet characters of the 
same modulus $ q \ (\geq 2) $. We define
\begin{eqnarray}
	&&\widehat{L}_{MT,j,r}(s_1,\cdots,s_j;s_{j+1},\cdots,s_{r+1};\chi_1,\cdots,\chi_r)
        \nonumber \\
        &&\qquad =  \sum_{m_1=1}^\infty \cdots \sum_{m_r=1}^\infty
        	\frac{\chi_1(m_1) \cdots \chi_r(m_r)}
                	{m_1^{s_1} \cdots m_j^{s_j}
                (m_1 + \cdots + m_j)^{s_{j+1}} 
                \cdots (m_1+ \cdots +m_r)^{s_{r+1}}}
     \label{eq:Miyagawa2} 
\end{eqnarray}
where $ 1 \leq j \leq r $ and $ s_1, \cdots s_r ,s_{r+1} $ are complex variables.
The series (\ref{eq:Miyagawa2}) is absolutely convergent in $ R_{j,r} $, 
and so $ \widehat{L}_{MT,j,r} $ is a regular function on $ R_{j,r} $.
\end{definition}
Definition \ref{def:L_MT,j,r} gives a generalized class which includes both
\begin{eqnarray}\label{eq:EZ-LL}
	\lefteqn{
                \cL_{EZ,r}(s_1,\cdots,s_r ;  \chi_1,\cdots,\chi_r)
                }	\nonumber \\
                && \qquad
                = \sum_{m_1=1}^\infty \cdots \sum_{m_r=1}^\infty
		\frac{\chi_1(m_1) \cdots \chi_r(m_r)}
                {m_1^{s_1} (m_1 + m_2)^{s_2} 
                \cdots (m_1 + \cdots + m_r)^{s_r}
        },
\end{eqnarray}
and
\begin{equation}\label{MT-L}
       	L_{MT,r}(s_1,\cdots,s_r;s_{r+1};  \chi_1,\cdots,\chi_r)
            =  \sum_{m_1=1}^\infty \cdots \sum_{m_r=1}^\infty
                \frac{\chi_1(m_1) \cdots \chi_r(m_r)}
                {m_1^{s_1} \cdots m_r^{s_r}
                (m_1 + \cdots + m_r)^{s_{r+1}}}
\end{equation}
as special cases. The series (\ref{eq:EZ-LL}) is introduced by Kamano {\cite{Kam}}, 
and he proved the meromorphic continuation to $ \C^r $. Also (\ref{MT-L}) is
introduced by Wu {\cite{Wu}} and he proved some analytic properties 
(see Theorem 3 in Matsumoto {\cite{Mat4}}).

\begin{remark}
Analytic properties of Apostol-Vu type (\ref{AV-zeta}) was also proved
by Okamoto \cite{Oka}, whose method is different from the method of 
Matsumoto \cite{Mat2} through the function (\ref{eq:Matsumoto2}). Okamoto's
method is based on the observation that (\ref{AV-zeta}) has the recursive structure
\begin{equation}\label{eq:induction} 
	{\zeta}_{AV,r} \longrightarrow
        {\zeta}_{AV,r-1} \longrightarrow 
        {\zeta}_{AV,r-2} \longrightarrow 
        \cdots \longrightarrow 
        {\zeta}_{AV,2} \longrightarrow 
        \zeta,
\end{equation}
where the right-most $ \zeta $ denotes the Riemann zeta-function.
Thus, analytic properties of (\ref{AV-zeta}) can be proved without using the
function (\ref{eq:Matsumoto2}) and the recursive structure (\ref{AV-RS}). 
\end{remark}

\begin{remark}
Matsumoto and Tanigawa {\cite{MaTa}} defined the multiple Dirichlet series 
\[
	\sum_{m_1=1}^\infty \cdots \sum_{m_r=1}^\infty
		\frac{a_1(m_1) a_2(m_2) \cdots a_r(m_r)}
                {m_1^{s_1} (m_1 + m_2)^{s_2} 
                \cdots (m_1 + \cdots + m_r)^{s_r}
        }
\]
which is a further generalization of (\ref{eq:EZ-LL}). They proved its 
several analytic properties. 
\end{remark}

\section{Statement of results}

\begin{theorem}\label{th:Main_Theorem1}
For $ 1 \leq j \leq r $, we have
\begin{enumerate}[$(i)$]
 \item the function $ \widehat{\zeta}_{MT,j,r}
 	(s_1,\cdots ,s_j; s_{j+1}, \cdots ,s_{r+1}) $ can be continued meromorphically 
	to the whole $ \C^{r+1} $-space,
 \item 
	in the case of $ j = r-1 $, the possible singularities of 
	$ \widehat{\zeta}_{MT,r-1,r} $
	are located only on the subsets of $ \C^{r+1} $ defined by one of the 
	following equations;
	\begin{eqnarray*}
		 && s_{r+1} = 1,	\\
         && s_j + s_r + s_{r+1} = 1 - \ell  
         	\quad (1 \leq j \leq r-1,\ \ell \geq -1),	\\
         && s_{j_1} + s_{j_2} + s_r + s_{r+1} = 2 - \ell 
         	\quad (1 \leq j_1 < j_2 \leq r-1,\ \ell \geq -1),	\\
         && \quad \vdots	\\
         && s_{j_1} + \cdots + s_{j_{r-2}} + s_r + s_{r+1} = r - 2 - \ell 
         	\quad (1 \leq j_1 < \cdots <j_{r-2} \leq r-1,
                				\ \ell \geq -1),	\\
         && s_1 + \cdots + s_{r-1} + s_r + s_{r+1} = r - 1 - d 
			\quad( d = -1,0,1,3,5,7,9, \cdots ). 
	 \end{eqnarray*}
	Also, in the cases of $ 1 \leq j \leq r-2 $, possible singularities of 
	$ \widehat{\zeta}_{MT,j,r} $ are located only 
	on the subsets of $ \C^{r+1} $ defined by one of the following equations;
  	\begin{eqnarray*}
		 && s_{r+1} = 1,	\\
         && s_r + s_{r+1} = 1 - d  \quad (d = -1,0,1,3,5,7,9, \cdots), \\
         && s_{r-1} + s_r + s_{r+1} = 3 - \ell \quad (\ell \in \N_0), \\
         && s_{r-2} + s_{r-1} + s_r + s_{r+1} = 4 - \ell 
         					\quad (\ell \in \N_0), \\
         && \qquad \vdots \\
         && s_{j+2} + s_{j+3} + \cdots + s_r + s_{r+1} = r-j-\ell
         					\quad (\ell \in \N_0), \\
         && s_{k_1} + s_{j+1} + \cdots + s_r + s_{r+1} = 1 - \ell'  
         	\quad (1 \leq k_1 \leq j,\ \ell' \geqq -(r-j)),	\\
         && s_{k_1} + s_{k_2} + s_{j+1} + \cdots +  s_r + s_{r+1} = 2 - \ell'
         	\quad (1 \leq k_1 < k_2 \leq j,\ \ell' \geqq -(r-j)), \\
         && \qquad \vdots	\\
         && s_{k_1} + \cdots + s_{k_{j-1}} + s_{j+1} + \cdots +  s_r + s_{r+1} 
         						= j - 1 - \ell' \\
         && \qquad \qquad \qquad \qquad \qquad \qquad \qquad
         		\quad (1 \leq k_1 < \cdots <k_{j-1} \leq j,
                				\ \ell' \geqq -(r-j)),    \\
         && s_1 + \cdots + s_j + s_{j+1} + \cdots +  s_r + s_{r+1} 
         						= j - \ell'
         	\quad (\ell' \geqq -(r-j)).
	 \end{eqnarray*}
 \item each of these singularities can be canceled by the corresponding linear
	factor, and
 \item $ \widehat{\zeta}_{MT,j,r} $ is of polynomial order with respect to
 	$ |\mathrm{Im}(s_{r+1})| $ .
\end{enumerate}
\end{theorem}

\medskip

\begin{theorem}\label{th:Main_Theorem2}
For $ 1 \leq j \leq r $, we have
\begin{enumerate}[$(i)$]
 \item the function 
	$ \widehat{L}_{MT,j,r}(s_1,\cdots ,s_j; s_{j+1}, 
		\cdots ,s_{r+1};\chi_1,\cdots \chi_r) $ 
	can be continued meromorphically to the $ \C^{r+1} $-space.
 \item If none of the characters $ \chi_1,\cdots,\chi_r $ are principal, 
	then $ \widehat{L}_{MT,j,r} $ is entire. If 
	$ \chi_{t_1},\cdots,\chi_{t_k} $ $ (1 \leq t_1 < \cdots < t_k \leq j) $ 
	and 
	$ \chi_{r-d_1},\cdots,\chi_{r-d_h} $ $ (1 \leq d_1 < \cdots < d_h \leq r-j) $
	are principal character and other characters are non-principal, 
	in the case of $ j = r-1 $, then possible singularities are located 
	only on the subsets of $ \C^{r+1} $ defined by one of the following equation;
    \begin{eqnarray}\label{possi.L.r-1}
         && s_{t_{u(1)}} + s_r + s_{r+1} = 1 - \ell  
         	\quad (1 \leq u(1) \leq k,\ \ell \geq - \delta_r),	\nonumber \\
         && s_{t_{u(1)}} + s_{t_{u(2)}} + s_r + s_{r+1} = 2 - \ell
         	\quad (1 \leq u(1) < u(2) \leq k,\ \ell \geq - \delta_r), \nonumber \\
         && \qquad \vdots	 \\
         && s_{t_{u(1)}} + \cdots + s_{t_{u(k-1)}} + s_r + s_{r+1} 
         						= k - 1 - \ell  \nonumber  \\
         && \qquad \qquad \qquad \qquad \qquad \qquad \qquad
         		\quad (1 \leq u(1) < \cdots < u(k-1) \leq k, \ \ell \geq - \delta_r),
			\nonumber \\
         && s_{t_1} + \cdots + s_{t_k} + s_r + s_{r+1} = k - \ell
         	\quad (\ell \geqq - \delta_r), \nonumber
	 \end{eqnarray}
	where
	\[
		\delta_r = \begin{cases}
					1 & (\chi_r \ is\  principal)	\\
					0 & (\chi_r \ is\  non\ principal)
					\end{cases},
	\]
	also in the cases of $ 1 \leq j \leq r-2 $, 
	then possible singularities are located only on the subsets of $ \C^{r+1} $
    	defined by one of the following equation;
	\begin{eqnarray}\label{possi.L.j}
        && s_{r-d_1 +1} + s_{r-d_1 +2} + \cdots + s_{r+1} = d_1 +1-\ell_0  
         					\quad (\ell_0 \in \N_0), \nonumber \\
        && \qquad \vdots \nonumber \\
        && s_{r-d_h +1} + s_{r-d_h +2} + \cdots + s_{r+1} = d_h +1-\ell_0
         					\quad (\ell \in \N_0), \nonumber \\
        && s_{t_{u(1)}} + s_{j+1} + \cdots + s_r + s_{r+1} = 1 - \ell'  
         	\quad (1 \leq u(1) \leq k,\ \ell' \geq - \Delta_j),	\nonumber \\
        && s_{t_{u(1)}} + s_{t_{u(2)}} + s_{j+1} + \cdots +  s_r + s_{r+1} = 2 - \ell' 
			\nonumber \\
        && \qquad \qquad \qquad \qquad \qquad \qquad \qquad
					(1 \leq u(1) < u(2) \leq k, \ \ell' \geq - \Delta_j),  \\
        && \qquad \vdots	\nonumber \\
        && s_{u(1)} + \cdots + s_{u(j-1)} + s_{j+1} + \cdots +  s_r + s_{r+1} 
         						= j - 1 - \ell' \nonumber \\
        && \qquad \qquad \qquad \qquad \qquad
         		\quad (1 \leq u(1) < \cdots <u(j-1) \leq k,
                				\ \ell' \geq - \Delta_j),	\nonumber \\
        && s_1 + \cdots + s_j + s_{j+1} + \cdots +  s_r + s_{r+1} 
         						= j - \ell'
         	\quad (\ell' \geq - \Delta_j), \nonumber
	 \end{eqnarray}
	 where $ \Delta_j = \delta_r + \delta_{r-1} + \cdots + \delta_{r-j} $.
	 Moreover, if $ \chi_r $ is principal character, then 
	 \[
		s_{r+1} = 1
	 \]
     is a possible singularity in addition to the above possible singularities
	 (\ref{possi.L.r-1}) and (\ref{possi.L.j}). 
 \item each of these singularities can be canceled by the corresponding linear
	factor, and
 \item $ \widehat{L}_{MT,j,r} $ is of polynomial order with respect to 
    $ |\mathrm{Im}(s_{r+1})| $.
\end{enumerate}
\end{theorem}

\medskip

\begin{remark}
In both Theorem \ref{th:Main_Theorem1} and Theorem \ref{th:Main_Theorem2}, 
the case $ j = r $ is known (see Theorem \ref{th:MT-zeta} and Theorem \ref{th:Wu1} below). 
It is interesting that the feature of possible singularities in the case $ j = r-1 $
is different from that in the cases $ 1 \leq j \leq r-2 $.
\end{remark}

\section{Proof of Theorem 1}

The proof of Theorem \ref{th:Main_Theorem1} and Theorem \ref{th:Main_Theorem2}
is similar to the argument of Matsumoto \cite{Mat1}, \cite{Mat2}, \cite{Mat}, 
\cite{Mat3}, \cite{Mat4}. The basic point is the use of the following integral 
representation.

\begin{lemma}\label{MBhyouji}
 We have 
 \begin{eqnarray}
        \lefteqn{\widehat{\zeta}_{MT,j,r}
                (s_1,\cdots s_j; s_{j+1},\cdots,s_r,s_{r+1})}	\nonumber \\
        && = \frac{1}{2\pi i}
        	\int_{(c)}
        	\frac{\Gamma(s_{r+1} + z) \Gamma(-z)}{\Gamma(s_{r+1})}		\nonumber \\
        && \label{zeta:MBhyouji} \qquad \qquad 
                \times \widehat{\zeta}_{MT,j,r-1}(s_1,\cdots,s_j;
                s_{j+1}, \cdots ,s_{r-1},s_r + s_{r+1} + z) \zeta(-z)dz 
                \label{M-B}
  \end{eqnarray}
 and
  \begin{eqnarray}
\lefteqn{\widehat{L}_{MT,j,r}(s_1,\cdots,s_j;s_{j+1},\cdots,s_{r+1};
				\chi_1, \cdots , \chi_r )}
        \nonumber \\
        &&= \frac{1}{2\pi i}
        	\int_{(c)}
        	\frac{\Gamma(s_{r+1} + z) \Gamma(-z)}{\Gamma(s_{r+1})}		\nonumber \\
        &&  \label{L:MBhyouji} \qquad
                \times \widehat{L}_{MT,j,r-1}(s_1,\cdots,s_j;
                s_{j+1}, \cdots, s_{r-1},s_r + s_{r+1} + z;
				\chi_1, \cdots , \chi_{r-1}) L(-z,\chi_r)dz,
\end{eqnarray}
where $ L(-z,\chi_r) $ is the Dirichlet $ L $-function attached to $ \chi_r $,
$ 1 \leq j \leq r-1 $ and $ - \mathrm{Re}(s_{r+1}) < c < -1 $ .
\end{lemma}
\textbf{Proof of Lemma \ref{MBhyouji}}. 
We prove only for $ \widehat{L}_{MT,j,r} $.
Using the Mellin-Barnes integral formula (\ref{eq:Mellin-Barnes}) for 
the multiple sum (\ref{eq:Miyagawa2}) 
with $ \lambda = m_r/(m_1 + \cdots + m_{r-1}) $ , 
we can formally obtain
\begin{eqnarray*}
\lefteqn{\widehat{L}_{MT,j,r}(s_1,\cdots,s_j;s_{j+1},\cdots,s_{r+1};
				\chi_1, \cdots , \chi_r)}
        \nonumber \\
        && =  \sum_{m_1=1}^\infty  \cdots
               	\sum_{m_r=1}^\infty
        	\frac{\chi_1(m_1) \cdots \chi_r(m_r)}
                {m_1^{s_1} \cdots m_j^{s_j}
                (m_1 + \cdots + m_j)^{s_{j+1}} 
                \cdots (m_1+ \cdots +m_r)^{s_{r}+s_{r+r}}}		\\
        && \qquad \qquad \qquad \qquad \qquad \qquad 
        	\qquad \qquad \qquad \qquad \qquad  \times
                \left( 1+ \frac{m_r}{m_1 + \cdots + m_{r-1}} \right)^{-s_{r+1}}	\\
        && = \sum_{m_1=1}^\infty  \cdots
               	\sum_{m_r=1}^\infty
        	\frac{\chi_1(m_1) \cdots \chi_r(m_r)}
                {m_1^{s_1} \cdots m_j^{s_j}
                (m_1 + \cdots + m_j)^{s_{j+1}} 
                \cdots (m_1+ \cdots +m_r)^{s_{r}+s_{r+r}}}			\\
	&& \qquad \qquad \qquad \qquad \qquad \qquad 
        	\times 
                \frac{1}{2\pi i}
        	\int_{(c)}
        	\frac{\Gamma(s_{r+1} + z) \Gamma(-z)}{\Gamma(s_{r+1})}
                \left( \frac{m_r}{m_1+ \cdots +m_{r-1}}  \right)^z dz	\\
	&&= \frac{1}{2\pi i}
        	\int_{(c)}
        	\frac{\Gamma(s_{r+1} + z) \Gamma(-z)}{\Gamma(s_{r+1})}	\\
	&&  \qquad 
                \times \sum_{m_1=1}^\infty \cdots \sum_{m_{r-1}=1}^\infty
        	\frac{\chi_1(m_1) \cdots \chi_{r-1}(m_{r-1})}{m_1^{s_1}\cdots m_j^{s_j}
                (m_1 + \cdots + m_j)^{s_{j+1}}
                \cdots (m_1 + \cdots + m_{r-1})^{s_r + s_{r+1}+z}
                }									\\
	&&  \qquad      \times \sum_{m_r = 1}^\infty \frac{\chi_r(m_r)}{{m_r}^{-z}} dz
\end{eqnarray*}
\begin{eqnarray*}
	&&= \frac{1}{2\pi i}
        	\int_{(c)}
        	\frac{\Gamma(s_{r+1} + z) \Gamma(-z)}{\Gamma(s_{r+1})}			\\
    &&  \qquad 
                \times \widehat{L}_{MT,j,r-1}(s_1,\cdots,s_j;
                s_{j+1}, \cdots, s_{r-1},s_r + s_{r+1} + z;
				\chi_1, \cdots , \chi_{r-1}) L(-z,\chi_r)dz.
\end{eqnarray*}
Now, we prove that $ \sum_{m=1}^\infty $ and $ \int_{(c)} $ can be exchanged. 
Put $ z = c + i w \ (-\infty < w < \infty) $. It is enough to prove that 
\begin{eqnarray*}
   I_{j,r} &=& \sum_{m_1=1}^\infty \cdots \sum_{m_r=1}^\infty
				\int_{-\infty}^{\infty}
				\left|
				\frac{\chi_1(m_1) \cdots \chi_r(m_r)}
                	{m_1^{s_1} \cdots m_j^{s_j}
                (m_1 + \cdots + m_j)^{s_{j+1}} 
                \cdots (m_1+ \cdots +m_{r-1})^{s_r + s_{r+1}}}	\right.	\\
			&& \qquad \qquad \qquad \qquad
				\left.  \times
				\left( \frac{m_r}{m_1 + \cdots + m_{r-1}} \right)^{c+iw}
                \frac{\Gamma(s_{r+1} + c + iw) \Gamma(- c - iw)}
                {\Gamma(s_{r+1})}
				\right| dw	\\   	
			&=& \widehat{\zeta}_{MT,j,r-1}(\sigma_1,\cdots, \sigma_j;
                \sigma_{j+1}, \cdots ,\sigma_{r-1} ,
                \sigma_r + \sigma_{r+1} + c) \zeta(-c)		\\
			&& \qquad \qquad \qquad \qquad
				\times \frac{1}{|\Gamma(s_{r+1})|}
                \int_{-\infty}^{\infty}
                |\Gamma(s_{r+1}+c+iw) \Gamma(-c-iw)|dw
\end{eqnarray*}
is bounded for each $ (s_1,s_2,\cdots,s_{r+1}) \in R_{j,r} $.
By using the Stirling's formula we have
\begin{eqnarray*}
\lefteqn{|\Gamma(s_{r+1}+c+iw) \Gamma(-c-iw)|} \\
 &=& \sqrt{2\pi} 
 	\left| 
 		\exp \left\{ \left( s_{r+1}+c+iw -\frac{1}{2} \right)
 		\log (s_{r+1}+c+iw) \right\}
 	\right|		\\
 & & \qquad \qquad \qquad \qquad
        \times 
        \left| 
 		\exp (-s_{r+1}-c-iw)
 	\right|  
        (1+O(|w|^{-1}))  \ \ \ (|w| \rightarrow \infty)\\
 &=& \sqrt{2\pi} \exp \{ -w \ \mathrm{arg}(s_{r+1}+c+iw) \}
 	\ O(|w|^{\sigma_{r+1} + c + \frac{1}{2}}) \\
 &=& O \left( \exp \left(-\frac{\pi}{2}|w| \right) \right),
\end{eqnarray*}
hence
\[
	\int_{-\infty}^{\infty} |\Gamma(s_{r+1}+c+iw) \Gamma(-c-iw)|dw
        = O(1) .
\]
This implies the assertion.
\qed

\medskip

These integral representations (\ref{zeta:MBhyouji}), (\ref{L:MBhyouji}), 
give the following inductive structure; 
\begin{equation}\label{eq:induction1} 
	\widehat{\zeta}_{MT,j,r} \longrightarrow 
        \widehat{\zeta}_{MT,j,r-1} \longrightarrow 
        \widehat{\zeta}_{MT,j,r-2} \longrightarrow 
        \cdots \longrightarrow 
        \widehat{\zeta}_{MT,j,j+1} \longrightarrow 
        \widehat{\zeta}_{MT,j,j} =
        \zeta_{MT,j},
\end{equation}
\begin{equation}\label{eq:induction2} 
	\widehat{L}_{MT,j,r} \longrightarrow 
        \widehat{L}_{MT,j,r-1} \longrightarrow 
        \widehat{L}_{MT,j,r-2} \longrightarrow 
        \cdots \longrightarrow 
        \widehat{L}_{MT,j,j+1} \longrightarrow 
        \widehat{L}_{MT,j,j} =
        L_{MT,j}.
\end{equation}

\bigskip

\begin{theorem}[K.Matsumoto \cite{Mat2}]\label{th:MT-zeta}
\begin{enumerate}[$(i)$]
\item The function $ \zeta_{MT,r}(s_1, \cdots, s_r; s_{r+1}) $ 
		can be meromorphically continued to the whole $ \C^{r+1} $-space.
\item The possible singularities of $ \zeta_{MT,r} $ are located only 
	  on the subsets of $ \C^{r+1} $ defined by one of the following equations;
	\begin{eqnarray*}
	 && s_j + s_{r+1} = 1 - \ell 
         	\quad (1 \leq j \leq r, \ \ell \in \N_0) ,		\\
 	 && s_{j_1} + s_{j_2} + s_{r+1} = 2-\ell 
         	\quad (1 \leq j_1 < j_2 \leq r, \ \ell \in \N_0) ,	\\
         && \qquad \ldots						\\
         && s_{j_1} + \cdots + s_{j_{r-1}} + s_{r+1} = r - 1 - \ell
         	\quad (1 \leq j_1 < \cdots <j_{r-1} \leq r, 
                \ \ell \in \N_0) ,					\\
         && s_1 + s_2 + \cdots + s_r + s_{r+1} = r ,
	\end{eqnarray*}
	where $ \N_0 $ denotes the set of non-negative integer.
\item Each of these singularities can be cancelled by the corresponding 
	  linear factor.

\item $ \zeta_{MT,r} $ is of polynomial order with 
	  respect to $ |\mathrm{Im}(s_{r+1})| $.
\end{enumerate}
\end{theorem}

\bigskip

\textbf{Proof of Theorem \ref{th:Main_Theorem1}}. 
When $ j = r $ the assertion is Theorem\ {\ref{th:MT-zeta}}. 
If $ j = r-1 $, (\ref{M-B}) implies
 \begin{eqnarray}
       		 && \widehat{\zeta}_{MT,r-1,r}
                (s_1,\cdots ,s_{r-1};s_r,s_{r+1}) \nonumber \\
       		 && = \frac{1}{2\pi i}
        	\int_{(c)}
        	\frac{\Gamma(s_{r+1} + z) \Gamma(-z)}{\Gamma(s_{r+1})}
                \zeta_{MT,r-1}(s_1,\cdots,s_{r-1};
                s_r + s_{r+1} + z) \zeta(-z) dz
                \label{eq:MBr-1,r}
 \end{eqnarray}
 where $ - \mathrm{Re}(s_{r+1}) < c < -1 $. By Theorem\ {\ref{th:MT-zeta}},
 the poles of $ \zeta_{MT,r-1}(s_1,\cdots,s_{r-1}; s_r + s_{r+1} + z) $
 as in a $z$-plane are
 \begin{eqnarray*}
	 && z = - s_j - s_r - s_{r+1} + 1 - \ell 
         	\qquad(1 \leq j \leq r-1, \ell \in \N_0), \\
 	 && z = - s_{j_1} - s_{j_2} - s_r - s_{r+1} + 2 - \ell 
         	\qquad(1 \leq j_1 < j_2 \leq r-1, \ell \in \N_0), \\
 	 && \qquad  \vdots \\
 	 && z = - s_{j_1} - \cdots - s_{j_{r-2}} - s_r - s_{r+1} + r-2 - \ell \\
         && \qquad \qquad \qquad \qquad \qquad \qquad \qquad \qquad \qquad 
         (1 \leq j_1 < \cdots < j_{r-1} \leq r-1, \ell \in \N_0), \\
	 && z = - s_1 - \cdots - s_{r-1} - s_r - s_{r+1} + r-1,
 \end{eqnarray*}
 all of which are located to the left of $ \mathrm{Re}(z) = c $.
 The other poles of the integrand on the right-hand side of (\ref{eq:MBr-1,r})
 are $ z = -s_{r+1}-n \ (n \in \N_0) $, $ z = n \ (n \in \N_0) $
 and $ z = -1 $.
 We shift the path of integration to the right to 
 $ \mathrm{Re}(z) = N - \varepsilon $, where $ N $ is a positive integer.
 Because $ \zeta_{MT,r-1}(s_1,\cdots,s_{r-1}; s_r) $ is of polynomial order 
 with respect to $ |\mathrm{Im}(s_r)| $, using Stirling's formula we obtain
 \begin{eqnarray*}
	&& \left|  \int_{c \pm iT}^{N-\varepsilon \pm iT}
        	\frac{\Gamma(s_{r+1} + z) \Gamma(-z)}{\Gamma(s_{r+1})}
                \zeta_{MT,r-1}(s_1,\cdots,s_{r-1};
                s_r + s_{r+1} + z) \zeta(-z) dz \right|	\\
	&&  \ll g(T) e^{-\pi T}	 \qquad (T \rightarrow \infty),
 \end{eqnarray*}
 where $ g $ is a certain polynomial. Hence, the shift of the path of integration
 is possible, and we obtain
 \begin{eqnarray} 
	\lefteqn{\widehat{\zeta}_{MT,r-1,r}
                (s_1,\cdots ,s_{r-1};s_r,s_{r+1})} \nonumber\\
	 && =\frac{1}{s_{r+1}-1}
            \zeta_{MT,r-1}
               (s_1,\cdots,s_{r-1};s_r + s_{r+1} - 1)	\nonumber\\
     &&	\quad	- \frac{1}{2}\zeta_{MT,r-1}
               (s_1,\cdots,s_{r-1};s_r + s_{r+1})  \nonumber\\
 	 && \quad + \sum_{n=1}^{\left[\frac{N}{2}\right]}
         	\binom{-s_{r+1}}{2n-1}
                \zeta_{MT,r-1}(s_1,\cdots,s_{r-1};s_r + s_{r+1} + 2n-1)
               \zeta(1-2n) \label{eq:MTr-1,r}			\nonumber\\ 
	 && \quad + 
        	\frac{1}{2\pi i}
        	\int_{(N - \varepsilon)}
        	\frac{\Gamma(s_{r+1} + z) \Gamma(-z)}{\Gamma(s_{r+1})}
                \zeta_{MT,r-1}(s_1,\cdots,s_{r-1};
                s_r + s_{r+1} + z) \zeta(-z) dz. 	
		\label{eq:shift zeta_{MT,r,r-1}}
 \end{eqnarray}
 The poles of the integrand of the last integral term is listed above, and
 hence we see that this integral is holomorphic at any points satisfying 
 all of the following inequalities;
 \begin{eqnarray*}
	 && \mathrm{Re}(s_{r+1}) > - N + \varepsilon, \nonumber\\
 	 && \mathrm{Re}(s_j + s_r + s_{r+1}) > 1 - N + \varepsilon, \nonumber\\
	 && \mathrm{Re}(s_{j_1} + s_{j_2} + s_r + s_{r+1})>2 - N + \varepsilon
         	\quad(1 \leq j_1 < j_2 \leq r-1), \nonumber	\\
 	 &&  \qquad\vdots 	\label{Region:holomor}\\
 	 && \mathrm{Re}(s_{j_1} + \cdots + s_{j_{r-2}} + s_r + s_{r+1})
         				> r - 2 - N + \varepsilon \nonumber \\
         && \qquad \qquad \qquad \qquad \qquad \qquad \qquad  \qquad \qquad
                	(1 \leq j_1 < \cdots < j_{r-2} \leq r-1), \nonumber \\
 	 && \mathrm{Re}(s_1 + \cdots + s_{r-1} + s_r + s_{r+1})
         				> r - 1 - N + \varepsilon. \nonumber
 \end{eqnarray*}
 Since $ N $ can be taken arbitrarily large, (\ref{eq:shift zeta_{MT,r,r-1}})
 implies the meromorphic continuation of 
 $ \widehat{\zeta}_{MT,r-1,r} (s_1,\cdots ,s_{r-1}; s_r,s_{r+1}) $
 to the whole $ \C^{r+1} $-space. 
 The first, the second and the third terms on right-hand side of 
 (\ref{eq:shift zeta_{MT,r,r-1}}) have a possible singularities that are 
 located only on the subsets of $ \C^{r+1} $ defined by one of the following 
 equations; 
 \begin{eqnarray*}
     && s_j + s_r + s_{r+1} + d = 1 - \ell  
         	\quad (1 \leq j \leq r-1,\ \ \ell \geqq 0),	\\
     && s_{j_1} + s_{j_2} + s_r + s_{r+1} + d = 2 - \ell 
         	\quad (1 \leq j_1 < j_2 \leq r-1,\ \ \ell \geqq 0),	\\
     && \quad \vdots	\\
     && s_{j_1} + \cdots + s_{j_{r-2}} + s_r + s_{r+1} + d = r - 2 - \ell 
       	\quad (1 \leq j_1 < \cdots <j_{r-2} \leq r-1,
                				\ \ \ell \geqq 0),	\\
     && s_1 + \cdots + s_{r-1} + s_r + s_{r+1} + d = r -1,
 \end{eqnarray*}
 where $ d = -1,0,1,3,5,7,\cdots \ (-1\leq d \leq N-1 ) $.
 Here, we note that
	$ \{ \ell + d \ | \ \ell \in \N_0, \ d = -1, 0, 1, 3, 5, \cdots \} 
		= \{ \ell \in \Z \ | \ \ell \geq -1 \} $.
 Since $ N $ can be arbitrarily large, we obtain the result in the case of
 $ j=r-1 $ in (ii).
 
 When $ j = r-2 $ in (\ref{zeta:MBhyouji}), and we shift the path of integration to
 the right to $ \mathrm{Re}(z) = N - \varepsilon $ to obtain 
 	\begin{eqnarray}
     \lefteqn{ \widehat{\zeta}_{MT,r-2,r}
                (s_1,\cdots ,s_{r-2};s_{r-1},s_r,s_{r+1})} \nonumber \\
  	 && \quad =
           \frac{1}{s_{r+1}-1}
            \widehat{\zeta}_{MT,r-1,r}
            (s_1,\cdots,s_{r-2};s_{r-1}, s_r + s_{r+1} - 1)	\nonumber \\
	 && \qquad \qquad \qquad \qquad
			- \frac{1}{2}\widehat{\zeta}_{MT,r-1,r}
            (s_1,\cdots,s_{r-2};s_{r-1}, s_r + s_{r+1})	\nonumber \\
 	 && \quad \quad + \sum_{n=1}^{\left[\frac{N}{2}\right]}
         	\binom{-s_{r+1}}{2n-1}
                \widehat{\zeta}_{MT,r-1,r}(s_1,\cdots,s_{r-2};s_{r-1},
                s_r + s_{r+1} + n)
               \zeta(1-2n) 		\label{eq:MTr-2,r} \\
	 && \quad \quad + 
        	\frac{1}{2\pi i}
        	\int_{(N - \varepsilon)}
        	\frac{\Gamma(s_{r+1} + z) \Gamma(-z)}{\Gamma(s_{r+1})}  \nonumber\\       
	 && \qquad \qquad \qquad \qquad 
			\times \widehat{\zeta}_{MT,r-1,r}(s_1,\cdots,s_{r-2};s_{r-1},
                s_r + s_{r+1} + z) \zeta(-z) dz.	\nonumber 
	\end{eqnarray}
   The possible singularities on the right-hand side of (\ref{eq:MTr-2,r}) are
        \begin{eqnarray*}
		 && s_{r+1} = 1,	\\
         && s_r + s_{r+1} + n = 1,  \\
         && s_j + s_{r-1} + s_r + s_{r+1} + n = 1 - \ell  
         	\quad (1 \leq j \leq r-2,\ \ \ell \geq -1),	\\
         && s_{j_1} + s_{j_2} + s_{r-1} + s_r + s_{r+1} + n = 2 - \ell 
         	\quad (1 \leq j_1 < j_2 \leq r-2,\ \ \ell \geq -1),	\\
         && \quad \vdots	\\
         && s_{j_1} + \cdots + s_{j_{r-3}} + s_{r-1} + s_r + s_{r+1} + n
                					 = r - 3 - \ell \\
         && \qquad \qquad \qquad \qquad \qquad \qquad \qquad
         	\quad (1 \leq j_1 < \cdots <j_{r-3} \leq r-2,
                				\ \ \ell \geq -1),	\\
         && s_1 + \cdots + s_{r-2} + s_{r-1} + s_r + s_{r+1} + n = r-2-d,
	\end{eqnarray*}
    where $ n,d = -1,0,1,3,5,7, \cdots \ (-1 \leq n \leq N) $. Since
    \begin{eqnarray*}
         && \{\ell + d  \ | \ \ell \in \{ -1 \} \cup \N_0,
        		\ d = -1,0,1,3,5, \cdots \} 
					= \{ \ell \in \Z \ |\  \ell \geq -2 \}, \\ 
         && \{ d + n  \ | \ d,n = -1,0,1,3,5, \cdots \ (-1 \leq n \leq N) \}
        		= \{ \ell \in \Z \ |\ \ell \geq -2 \},
    \end{eqnarray*}
    the above possible singularities can be rewritten as follows;
    \begin{eqnarray*}
		 && s_{r+1} = 1,	\\
         && s_r + s_{r+1} = 1 - n, \\
         && s_j + s_{r-1} + s_r + s_{r+1} = 1 - \ell  
         	\quad (1 \leq j \leq r-2,\ \ \ell \geqq -2),	\\
         && s_{j_1} + s_{j_2} + s_{r-1} + s_r + s_{r+1} = 2 - \ell 
         	\quad (1 \leq j_1 < j_2 \leq r-2,\ \ \ell \geqq -2),	\\
         && \quad \vdots	\\
         && s_{j_1} + \cdots + s_{j_{r-3}} + s_{r-1} + s_r + s_{r+1} 
         	= r - 3 - \ell \\
         && \qquad \qquad \qquad \qquad \qquad \qquad \qquad
         	\quad (1 \leq j_1 < \cdots <j_{r-2} \leq 
                				r-2,\ \ \ell \geqq -2),	\\
         && s_1 + \cdots + s_{r-2} + s_{r-1} + s_r + s_{r+1} = r - 2 - \ell 
         					\quad (\ell \geqq -2).
	 \end{eqnarray*}
     Since $ N $ can be taken arbitrarily large, we obtain the results of (ii)
	 in the case of $ j = r-2 $.
	 

Let $ k = r-j \ (k \geq 2) $. Assume that the assertion of 
Theorem \ref{th:Main_Theorem1} is true in the case of $ r-j = 2,3, \cdots ,k-1 $, 
and we prove by induction the assertion in the case of $ r-j = 2 $.
By Lemma \ref{MBhyouji}, we obtain
\begin{eqnarray}
        	&& \widehat{\zeta}_{MT,r-k,r}
                	(s_1,\cdots ,s_{r-k}; s_{r-k+1},\cdots,s_r,s_{r+1}) \nonumber\\
        	&& = \frac{1}{2\pi i}
        		\int_{(c)}
        		\frac{\Gamma(s_{r+1} + z) \Gamma(-z)}{\Gamma(s_{r+1})}
                       \nonumber \\
                && \qquad \qquad
 			\times
                        \widehat{\zeta}_{MT,r-k,r-1}(s_1,\cdots,s_{r-k};
		s_{r-k+1}, \cdots ,s_{r-1}, s_r + s_{r+1} + z) \zeta(-z)dz,
		\label{MB-MT,r-k,r}
\end{eqnarray}
where $ 1 \leq j \leq k-1, \  - \mathrm{Re}(s_{r+1}) < c < -1 $.
By assumption of induction, we find that the possible singularities of 
 $ \widehat{\zeta}_{MT,r-k,r-1}(s_1,\cdots,s_{r-k};
		s_{r-k+1}, \cdots ,s_{r-1} , s_r + s_{r+1} + z) $
as a function in $ z $ are
\begin{eqnarray*}
	 && z = -s_r -s_{r+1} +1, \\
	 && z = -s_{r-1}-s_r-s_{r+1}+2-d \quad (d=-1,0,1,3,5,7,\cdots), \\
	 && z = -s_{r-2}-s_{r-1}-s_r-s_{r+1}+3-\ell \quad(\ell \in \N_0),\\
	 && \qquad \vdots \\
	 && z = -s_{r-k+2}-\cdots-s_{r-1}-s_r-s_{r+1}+k-1-\ell 
         					\quad(\ell \in \N_0), \\
	 && z = -s_{j_1}-s_{r-k+1}-\cdots-s_{r-1}-s_r-s_{r+1}+1-\ell'\\
         && \qquad \qquad \qquad \qquad \qquad \qquad
         		\qquad(1\leq j_1 \leq r-k,\  \ell'\geqq -k+2), \\
	 && z = -s_{j_1}-s_{j_2}-s_{r-k+1}-\cdots-s_{r-1}
         					-s_r-s_{r+1}+2-\ell'\\
         && \qquad \qquad \qquad \qquad \qquad \qquad
         	\quad(1\leq j_1 < j_2 \leq r-k,\  \ell'\geqq -k+2), \\
	 && \qquad \vdots \\
	 && z = -s_{j_1}-\cdots-s_{j_{r-k-1}}-s_{r-k+1}-\cdots-s_{r-1}
         					-s_r-s_{r+1}+r-k-\ell'\\
         && \qquad \qquad \qquad \qquad \qquad \qquad
         	(1\leq j_1 < \cdots < j_{r-k-1} \leq r-k, \ell'\geqq -k+2), \\
	 && z = -s_1-s_2-\cdots-s_{r-k+1}-\cdots-s_{r-1}
         					-s_r-s_{r+1}+r-k+1-\ell' \\
         && \qquad \qquad \qquad \qquad \qquad \qquad
         	\qquad \qquad \qquad \qquad \qquad \qquad 
                \qquad(\ell' \geqq -k+2), 
\end{eqnarray*}
all of which are located to the left of $ \mathrm{Re}(z) = c $.
The other poles of the integrand on the right-hand side of 
(\ref{MB-MT,r-k,r}) are $ z = -s_{r+1} - n \ (n \in \N_0) $
, $ z = n \ (n \in \N_0) $ and $ z = -1 $.
We shift the path of integration to the right to 
$ \mathrm{Re}(z) = N - \varepsilon $, where $ N $ is a positive integer.
Since the shift of the path of integration is possible as before, we obtain
 \begin{eqnarray}
	\lefteqn{\widehat{\zeta}_{MT,r-k,r}
                (s_1,\cdots ,s_{r-k};s_{r-k+1},\cdots,s_{r+1})} \nonumber \\
         && \nonumber \\
	 &=&    \frac{1}{s_{r+1}-1}
              	\widehat{\zeta}_{MT,r-k,r-1}
                (s_1,\cdots ,s_{r-k};s_{r-k+1},\cdots
		,s_{r-1},s_r + s_{r+1} - 1) \nonumber \\
         && \qquad \qquad \qquad \qquad \qquad
               - \frac{1}{2} 
               \widehat{\zeta}_{MT,r-k,r-1}
                (s_1,\cdots ,s_{r-k};s_{r-k+1},\cdots
		,s_{r-1},s_r + s_{r+1})   \nonumber \\
 	 && \quad + \sum_{n=1}^{\left[\frac{N}{2}\right]}
         	\binom{-s_{r+1}}{2n-1}
        	\widehat{\zeta}_{MT,r-k,r-1}
                (s_1,\cdots ,s_{r-k};s_{r-k+1},\cdots
		,s_{r-1},s_r + s_{r+1} + 2n - 1)	\nonumber \\
         && \qquad \qquad \qquad \times 
         	\zeta(1-2n) \label{eq:j,r} 
                \\
	 && \quad + 
        	\frac{1}{2\pi i}
        	\int_{(N - \varepsilon)}
        	\frac{\Gamma(s_{r+1} + z) \Gamma(-z)}{\Gamma(s_{r+1})}	
                \nonumber \\
         && \qquad \qquad \qquad \times
                \widehat{\zeta}_{MT,r-k,r-1}
                (s_1,\cdots ,s_{r-k};s_{r-k+1},\cdots
		,s_{r-1},s_r + s_{r+1} + z) \zeta(-z) dz.	\nonumber
	\end{eqnarray}
The first, the second and the third terms on right-hand side of (\ref{eq:j,r}) 
have a possible singularities that are located only on the subsets of 
$ \C^{r+1} $ defined by one of the following equations;
\begin{eqnarray}
	 &&   
		  s_{r+1} = 1,	\nonumber\\
     && 
          s_r + s_{r+1} + n = 1,	\nonumber\\
     && 
          s_{r-1} + s_r + s_{r+1} + n = 2 - d 
          	\quad (d = -1,0,1,3,5,7, \cdots), \nonumber\\
     && 
          s_{r-2} + s_{r-1} + s_r + s_{r+1} + n = 3 - \ell 
          	\quad (\ell \in \N_0), \nonumber\\
     &&  
         	\qquad \vdots \nonumber\\
     && 
          s_{r-k+2} + s_{r-k+3} + \cdots + s_r + s_{r+1} +n = k-1-\ell
         				\quad (\ell \in \N_0), \nonumber\\
     && 
          s_{j_1} + s_{r-k+1} + \cdots + s_r + s_{r+1}+n = 1 - \ell' \\
     && 
         	\qquad \qquad \qquad \qquad \qquad \qquad
         	\ \  (1 \leq j_1 \leq r-k,\ \ell' \geqq -(k-1)),	
                \label{eq:possible sing}\nonumber\\ 
     && 
          s_{j_1} + s_{j_2} + s_{r-k+1} + \cdots +  s_r + s_{r+1} +n
          	= 2 - \ell'	\nonumber\\
     && 
          \qquad \qquad \qquad \qquad \qquad \qquad
           	\quad (1 \leq j_1 < j_2 \leq r-k,\ \ \ell' \geqq -(k-1)), 
                \nonumber\\
     &&  
         	\qquad \vdots	\nonumber\\
     && 
          s_{j_1} + \cdots + s_{j_{r-k-1}} + s_{r-k+1} 
          		+ \cdots +  s_r + s_{r+1}+n = r-k-1- \ell' \nonumber\\
     && 
          \qquad \qquad \qquad \qquad \quad
         \quad (1 \leq j_1<\cdots<j_{r-k-1} \leq r-k,\ \ \ell' \geqq -(k-1)),
                                                         \nonumber\\
     && 
          s_1 + \cdots + s_{r-k} + s_{r-k+1} + \cdots +  s_r + s_{r+1} +n
         					= r-k - \ell'	\nonumber\\
     && 
          \qquad \qquad \qquad \qquad \qquad \qquad \qquad \qquad \qquad
          \qquad \qquad
         	\quad (\ell' \geqq -(k-1)) , \nonumber
\end{eqnarray}
where $ n = -1,0,1,3,5,7,\cdots \ (1 \leq n \leq N-1) $.
The last integral of (\ref{eq:j,r}) is holomorphic at any satisfying all of the 
following inequalities;
  \begin{eqnarray}
	 && \mathrm{Re}(s_{r+1}) > - N + \varepsilon,  \nonumber \\ 
 	 && \mathrm{Re}(s_r + s_{r+1}) > 1 - N + \varepsilon, \nonumber \\
     && \mathrm{Re}(s_{r-1} + s_r + s_{r+1}) > 2- N + \varepsilon,  \nonumber \\
     &&  \qquad \vdots \nonumber \\
	 && \mathrm{Re}(s_{r-k+2} + s_{r-k+3} + 
            			\cdots + s_r + s_{r+1}) > k-1- N + \varepsilon,
         		\nonumber \\
 	 && \mathrm{Re}(s_{j_1} + s_{r-k+1}+ \cdots + s_r + s_{r+1})
         					 > k- N + \varepsilon 
                 \quad (1 \leq j_1 \leq r-k),	\label{(4.8)} \\
     && \mathrm{Re}(s_{j_1} + s_{j_2} + s_{r-k+1}+ \cdots + s_r + s_{r+1})
         					 > k+1- N + \varepsilon \nonumber \\
     && \qquad \qquad \qquad \qquad \qquad \qquad \qquad \qquad \qquad \qquad
									(1 \leq j_1 < j_2 \leq r-k),	\nonumber \\
     &&  \qquad\vdots \nonumber \\
     && 
        \mathrm{Re}(s_{j_1} + \cdots + s_{j_{r-k-2}} 
        + s_{r-k+1}+ \cdots+ s_r + s_{r+1}) > r-2 - N+\varepsilon \nonumber \\
     && \qquad \qquad \qquad \qquad \qquad \qquad \qquad \qquad \qquad
         		\quad(1 \leq j_1 < \cdots < j_{r-k-2} \leq r-k), \nonumber \\
 	 && \mathrm{Re}(s_1 + \cdots + s_{r-k} + s_{r-k+1} 
         		+ \cdots + s_r + s_{r+1}) > r-1- N + \varepsilon. \nonumber 
	\end{eqnarray}
Since $ N $ can be taken arbitrarily large, (\ref{(4.8)}) implies the 
meromorphic continuation of 
$  \widehat{\zeta}_{MT,r-k,r} (s_1,\cdots ,s_{r-k}; s_{r-k+1}, \cdots ,s_{r+1}) $ 
to the whole $ \C^{r+1}$ space.
By the method similar to that as in the case of $ j = r - 2 $, 
we obtain the result in the case of $ 2 \leq j \leq r $ in (ii).

Let
 \begin{eqnarray*}
	\lefteqn{\Phi_{r-k,r,N}(s_1,\cdots,s_{r-k};s_{r-k+1},\cdots,s_{r+1})}\\
	 &=& (s_{r+1}-1) 
         \prod_{{-1 \leq d \leq N-1} \atop {d:0 \ or \ \mathrm{odd}}} 
         	(s_r + s_{r+1} - 2 + d) \\
	 && \times \prod_{\ell = 0}^{N}
         	\{ (s_{r-1} + s_r + s_{r+1} - 3 - \ell)
                   (s_{r-2} + s_{r-1} + s_r + s_{r+1} - 4 - \ell) \\
         && \qquad \times \cdots \times 
         	(s_{r-k+1} + \cdots + s_r + s_{r+1} -k-1+\ell)\} \\
	 && \times \prod_{\ell' = -k}^{N} \Bigg\{
        	\prod_{j_1=1}^{r-k}(s_{j_1}+s_{r-k+1} 
                	+ \cdots + s_r + s_{r+1} -1+\ell')	\\
         && \qquad \times
         	\prod_{1 \leq j_1<j_2 \leq r-k}^{r-k}(s_{j_1} + s_{j_2}
                + s_{r-k+1} + \cdots + s_r + s_{r+1} -1+\ell')	\\
	 && \qquad \times \cdots	\\
         && \qquad \times
         	\prod_{1 \leq j_1< \cdots <j_{r-k-1} \leq r-k}^{r-k}
                (s_{j_1} + \cdots +s_{j_{r-k-1}} 
                + s_{r-k+1} + \cdots + s_r + s_{r+1} -1+\ell')	\\
         && \qquad \times
         	(s_1 + \cdots + s_{r-k} + s_{r-k+1} 
                + \cdots + s_r + s_{r+1} -r + k + \ell')
         	\Bigg\}
\end{eqnarray*}
where $ N $ is positive integer. By (\ref{eq:j,r}) and (ii), 
        \[
        	\widehat{\zeta}_{MT,r-k,r}
                (s_1,\cdots ,s_{r-k};s_{r-k+1},\cdots,s_{r+1})
                \Phi_{r-k,r,N}(s_1,\cdots,s_{r-k};s_{r-k+1},\cdots,s_{r+1})
        \]
is shown to be holomorphic, to obtain (iii). Finally we can also prove
(iv) also by the induction assumption on the order 
$ \widehat{\zeta}_{MT,r-k,r-1} $ and Stirling's formula. 
Hence the proof of Theorem \ref{th:Main_Theorem1} is complete.
\qed

\section{Proof of Theorem 2}

\begin{theorem}[Wu \cite{Wu}] \label{th:Wu1}
 The function $ L_{MT,r}(s_1,\cdots,s_r ; s_{r+1} ; \chi_1,\cdots,\chi_r) $
 can be meromorphically continued to the $ \C^{r+1} $-space.
 If none of the characters $ \chi_1, \cdots, \chi_r $ are principal, then
 $ L_{MT,r} $ is entire.
 If there are $ k $ principal characters $ \chi_{t_1},\cdots,\chi_{t_k} $
 among them, then possible singularities are located only on the subsets
 of $ \C^{r+1} $ defined by one of the following equations;
 \begin{eqnarray*}
	 && s_{t_{u(1)}} + s_{r+1} = 1 - \ell 
         	\quad (1 \leq u(1) \leq k, \ \ell \in \N_0) ,		\\
 	 && s_{t_{u(1)}} + s_{t_{u(2)}} + s_{r+1} = 2-\ell 
         	\quad (1 \leq u(1) < u(2) \leq k, \ \ell \in \N_0) ,	\\
         && \qquad \vdots						\\
         && s_{t_{u(1)}} + \cdots + s_{t_{u(k-1)}} + s_{r+1} = k - 1 - \ell \\
         && \qquad \qquad \qquad \qquad \qquad \qquad
			 (1 \leq u(1) < \cdots < u(k-1) \leq k, \ \ell \in \N_0) ,	\\
         && s_{t_1} + s_{t_2} + \cdots + s_{t_k} + s_{r+1} 
				= k - \ell \left( 1 - \left[ \frac{k}{r} \right] \right) 
			\quad (\ell \in \N_0),
	\end{eqnarray*}
 where $ 1 \leq h \leq k, 1 \leq u(1) < \cdots < u(h) \leq k, \ell \in N_0 $. 
\end{theorem}


\textbf{Proof of Theorem \ref{th:Main_Theorem2}}. 
For (iii) and (iv) the method is exactly the same as in the proof of 
Theorem \ref{th:Main_Theorem1}.
When $ j = r $ the assertion is nothing but Theorem \ref{th:Wu1}. 
If $ j = r-1 $, (\ref{L:MBhyouji}) implies 
\begin{eqnarray}
\lefteqn{\widehat{L}_{MT,r-1,r}(s_1,\cdots,s_{r-1};s_r, s_{r+1};
				\chi_1, \cdots , \chi_r )}
        \nonumber \\
        &&= \frac{1}{2\pi i}
        	\int_{(c)}
        	\frac{\Gamma(s_{r+1} + z) \Gamma(-z)}{\Gamma(s_{r+1})}		\nonumber \\
        &&  \qquad \qquad \label{L_r-1,r:MBhyouji}  
                \times L_{MT,r-1}(s_1,\cdots, s_{r-1}; s_r + s_{r+1} + z;
				\chi_1, \cdots , \chi_{r-1}) L(-z,\chi_r)dz
\end{eqnarray}
where $ -\mathrm{Re}(s_{r+1}) < c < -1 $ , and $ L(¥,\chi_r) $ is 
Dirichlet $L$-function. By Theorem \ref{th:Wu1}, the poles of 
$ L_{MT,r-1}(s_1,\cdots, s_{r-1}; s_r + s_{r+1} + z; \chi_1, \cdots , \chi_{r-1}) $
as in the $ z$-plane are located to the left of $ \mathrm{Re}(z) = c $ .
The other poles of the integrand on the right-hand side of (\ref{L_r-1,r:MBhyouji}) 
are $ z = -s_{r+1} - n \ (n \in \N_0), \ z = n \ (n \in \N_0) $.
Also, when $ \chi_r $ is principal, $ z = -1 $ is a simple pole.
We shift the path of integration of (\ref{L_r-1,r:MBhyouji}) to the right to 
$ \mathrm{Re}(z) = N - \varepsilon $ , to obtain
\begin{eqnarray}
\lefteqn{\widehat{L}_{MT,r-1,r}(s_1,\cdots,s_{r-1};s_r, s_{r+1};
				\chi_1, \cdots , \chi_r )}
        \nonumber \\
        &&= \frac{1}{s_{r+1}-1} L_{MT,r-1}(s_1,\cdots, s_{r-1}; s_r + s_{r+1} - 1;
				\chi_1, \cdots , \chi_{r-1})\cdot 
				\frac{\varphi(q)}{q} \cdot \delta_r  \nonumber		\\
		&&	\qquad
			+ \sum_{n=0}^{N-1} \binom{-s_{r+1}}{n}
			L_{MT,r-1}(s_1,\cdots, s_{r-1}; s_r + s_{r+1} + n; \chi_1, \cdots , \chi_{r-1})
			L(-n, \chi_r)   \label{L_r-1,r:MBhyouji2}  \\
		&&	\qquad + \frac{1}{2\pi i}
        	\int_{(N - \varepsilon)}
        	\frac{\Gamma(s_{r+1} + z) \Gamma(-z)}{\Gamma(s_{r+1})}		\nonumber \\
        &&  \qquad \qquad  
                \times L_{MT,r-1}(s_1,\cdots, s_{r-1}; s_r + s_{r+1} + z;
				\chi_1, \cdots , \chi_{r-1}) L(-z,\chi_r)dz	\nonumber		
\end{eqnarray}
where $ \delta_r $ is defined in the statement of Theorem \ref{th:Main_Theorem2}.
Futher, if $ \chi_{t_1},\ \cdots ,\ \chi_{t_k} \ ( 1 \leq t_1 < \cdots < t_k \leq r-1 ) $ 
are principal and the others are non-principal, possible singularities of
(\ref{L_r-1,r:MBhyouji2}) are
\begin{eqnarray}
        && s_{t_{u(1)}} + s_r + s_{r+1} = 1 - \ell  
         	\quad (1 \leq u(1) \leq k,\ \ \ell \geq - \delta_r),	\nonumber \\
        && s_{t_{u(1)}} + s_{t_{u(2)}} + s_r + s_{r+1} = 2 - \ell 
         	\quad (1 \leq u(1) < u(2) \leq k,\ \ \ell \geq - \delta_r),	\nonumber \\
        && \quad \vdots		\label{L_r-1,r_possi.sing.}	\\
        && s_{t_{u(1)}} + \cdots + s_{t_{u(k-1)}} + s_r + s_{r+1} = k - 1 - \ell \nonumber \\
        && \qquad \qquad \qquad \qquad \qquad \qquad 
         	(1 \leq u(1) < \cdots <u(k-1) \leq k, \ \ \ell \geq - \delta_r),	\nonumber \\
        && s_{t_1} + \cdots + s_{t_k} + s_r + s_{r+1} = k - \ell
			\quad (\ell \geq - \delta_r), \nonumber
\end{eqnarray}
moreover, if $ \chi_{t_1},\ \cdots ,\ \chi_{t_k} \ ( 1 \leq t_1 < \cdots < t_k \leq r-1 ) $ 
and $ \chi_r $ are principal and the others are non-principal, then 
\[
	s_{r+1} = 1
\]
is also a possible singularity.
Proof in the case of $ 1 \leq j \leq r-2 $ is the same as the proof of 
Theorem \ref{th:Main_Theorem1}; we can prove the assertion using the induction on
$ k $ with $ k = r-j $. Also, how to deal with Dirichlet characters is similar to 
the case of $ j = r-1 $ .
\qed

\bigskip
\author{Takashi Miyagawa}:	\\
Graduate School of Mathematics, \\
Nagoya University, \\
Chikusa-ku, Nagoya, 464-8602 Japan	\\
E-mail: d15001n@math.nagoya-u.co.jp

\end{document}